\documentclass[11pt]{article}

\title{\vskip-1.0em Triviality of the generalized Lau product associated to a Banach algebra homomorphism}

\author{Yemon Choi}
\date{December 2, 2015}

\usepackage{amsfonts,amsthm}
\newcounter{pulse}
\theoremstyle{plain}

\newtheorem{prop}[pulse]{Proposition}

\theoremstyle{definition}
\newtheorem{dfn}[pulse]{Definition}
\newtheorem{eg}[pulse]{Example}
\newtheorem{rem}[pulse]{Remark}

\newcommand{\dt}[1]{{\itshape#1}}
\newcommand{\iso}{\cong}
\newcommand{\fu}[1]{{#1}^\sharp}

\begin{document}
\maketitle
\begin{abstract}
Several papers have, as their raison d'\^{e}tre, the exploration of the \dt{generalized Lau product} associated to a homomorphism $T:B\to A$ of Banach algebras. In this short note, we demonstrate that the generalized Lau product is isomorphic as a Banach algebra to the usual direct product $A\oplus B$. We also correct some misleading claims made about the relationship between this generalized Lau product, and an older construction of Monfared (Studia Mathematica, 2007).

\medskip
\noindent MSc 2010 classification: 46H20
\end{abstract}

\bigskip

\begin{section}{The generalized Lau product}
Motivated by a construction of M. S. Monfared \cite{Mon2007}, several authors have in recent years written papers on a certain construction, which manufactures a Banach algebra $A\times_T B$ given a pair of Banach algebras $(A,B)$ and a continuous algebra homomorphism $T:B\to A$.

\begin{dfn}\label{d:genlauprod}
The underlying Banach space of $A\times_T B$ is the usual product/sum $A\oplus B$; multiplication is defined by the following rule:
\[ (a_1,b_1) \bullet_T (a_2,b_2) = (a_1a_2 + T(b_1)a_2 + a_1T(b_2)\,, b_1b_2) \]
It is clear that $\bullet_T$ is continuous, and easily verified by hand that it is associative (although this will also follow from the proof of Proposition~\ref{p:kill it with fire}). By renorming if necessary to obtain a submultiplicative norm, one obtains a Banach algebra $A\times_T B$.
\end{dfn}

The construction in Definition~\ref{d:genlauprod} has gone by various names: ``generalized Lau product''; ``the $T$-Lau product''; or ``the Lau product defined by a Banach algebra morphism''. This last phrase is the one used in \cite{BhDa2013}, which appears to be the earliest occurrence of this construction. Since \cite{BhDa2013} there have been several papers on the theme of the ``generalized Lau product'', by various authors: see
\cite{MR3392255},
\cite{MR3294267},
%\cite{MR3040708},
\cite{MR3376847},
\cite{JaNe_pre},
\cite{JaNe2014},
\cite{MR3166398},
\cite{MR3287161},
\cite{MR3296113},
\cite{PoRa_pre2},
\cite{PoRa_pre}.

\begin{rem}\label{r:kingsley says no}
In \cite{BhDa2013} it is claimed that this construction extends that of~\cite{Mon2007}, and that it might have some bearing on extensions of $C^*$-algebras in the sense of BDF theory. The second claim is not given clear justfication in the article; and the first claim is misleading, to say the least, as will be explained in Section~\ref{s:not Monfared}).
\end{rem}

Unfortunately, the following elementary observation casts doubt on the whole enterprise. We leave the proof to the reader.

\begin{prop}\label{p:kill it with fire}
Let $A\oplus_{\rm alg} B$ denote the usual sum of Banach algebras, i.e.~we equip the Banach spaces $A\oplus B$ with co-ordinatewise product
\[ (a_1,b_1)\cdot (a_2,b_2) = (a_1a_2, b_1b_2). \]
Define $\phi: A\times_T B \to A\oplus_{\rm alg} B$ by
$\phi(a,b) = (a+T(b), b)$.
Then $\phi$ is continuous, linear and bijective; and $\phi((a,b)\bullet_T(c,d)) = \phi(a,b)\cdot\phi(c,d)$.
In particular, $A\times_T B$ is isomorphic as a Banach algebra to $A\oplus_{\rm alg} B$.
\end{prop}

\begin{rem}
The proof of Proposition~\ref{p:kill it with fire} was found when the present author was reading the preprint \cite{PoRa_pre}, in particular the proof of Proposition 4.1 in that paper. It seems that the authors of \cite{PoRa_pre} noticed something very similar to the isomorphism of Proposition~\ref{p:kill it with fire}, but only in a restricted setting. Similarly, it is observed just after Theorem 4.2 of \cite{JaNe_pre} that there is an isomorphism $A\times_T B\to A\oplus_{\rm alg} B$, but the authors of \cite{JaNe_pre} only state this in the case where $A$ and $B$ are commutative and semisimple, and do not mention that Proposition~\ref{p:kill it with fire} applies in full generality.
\end{rem}

 Hopefully, any future attempts to study the generalized Lau product, as defined in Definition~\ref{d:genlauprod}, will bear in mind that $A\times_T B \iso A\oplus_{\rm alg} B$ regardless of the choice of $T$.
It is a basic theme, when defining a property of Banach algebras, to see how it behaves under forming binary sums/products of algebras. For many of the properties considered in the items of the bibliography, stability of such properties under forming binary sums/products is either known or refuted by old work.

\end{section}

\begin{section}{Comparison with an older construction}\label{s:not Monfared}
It is claimed in \cite{BhDa2013} that the construction presented there extends the one studied by Monfared in~\cite{Mon2007}. If this were true, then Monfared's construction would be a special case of Definition~\ref{d:genlauprod}, and hence would be trivial for the same reason.

In fact, it is \emph{not} true that Monfared's construction is a special case of Definition~\ref{d:genlauprod}. In this section we shall briefly explain why.

\begin{dfn}[Monfared, \cite{Mon2007}]
\label{d:lauprod}
Let $A$ and $B$ be Banach algebras and let $\varphi:B\to {\mathbb C}$ be a character, i.e.~a non-zero homomorphism. 
The \dt{Lau product} of $A$ and $B$ with respect to $\varphi$ is defined to be the Banach algebra whose underlying Banach space is $A\oplus B$, equipped with the multiplication operation
\[ (a_1,b_1) \bullet_\varphi (a_2,b_2) = (a_1a_2 + \varphi(b_1)a_2 + a_1\varphi(b_2)\,, b_1b_2)\,. \]
\end{dfn}

As observed by Bhatt and Dabhi, it is immediate that when $A$ has an identity element $e_A$ and $\varphi$ is a character on $B$, we may define a continuous homomorphism $T :B\to A$ by $T(a)=\varphi(b)e_A$. If we do this, then the Lau product of $A$ and $B$ with respect to $\varphi$ does indeed coincide with $A\times_T B$, and hence by Proposition~\ref{p:kill it with fire} it is isomorphic to the usual direct product of algebras $A\oplus_{\rm alg} B$.

What seems to have gone unremarked in \cite{BhDa2013} is that when $A$ does not have an identity element, Definition~\ref{d:lauprod} is no longer a special case of Definition~\ref{d:genlauprod}.
Indeed, in general the Lau product of $A$ and $B$ is not isomorphic as an algebra to $A\oplus_{\rm alg} B$. This is not a new observation --- it is implicit in \cite{Mon2007} --- but for sake of completeness we shall give an example to illustrate this.

\begin{eg}\label{eg:different}
Let $A$ be an arbitrary Banach algebra. Take $B={\mathbb C}$ and let ${\rm id}$ denote the identity map on ${\mathbb C}$. It was observed in \cite{Mon2007} (and it is also clear from Definition~\ref{d:lauprod}) that the corresponding Lau product of $A$ and $B$ is isomorphic to $\fu{A}$, the usual unitization of~$A$. Observe that $A\oplus_{\rm alg} {\mathbb C}$ has an identity element if and only if $A$ does. Hence, if $A=c_0({\mathbb N})$ (for example), then $A\oplus_{\rm alg}{\mathbb C}$ cannot be isomorphic to $\fu{A}$.
\end{eg}

A minor variation on the proof of Proposition~\ref{p:kill it with fire} shows that given $A,B,\varphi$ as in Definition~\ref{d:lauprod}, there is a natural way to identify $A\times_\varphi B$ with a closed, codimension-one subalgebra (not an ideal, in general!) of $\fu{A}\oplus_{\rm alg} B$. Given that there continue to be papers exploring the Lau product (as in Definition~\ref{d:lauprod}), it therefore seems worthwhile to record the following consequence, whose proof we omit since it is straightforward.

\begin{prop}
Let $Q$ be a property of Banach algebras which is preserved by all isomorphisms of Banach algebras (not just the isometric ones). Suppose also that:
\begin{itemize}
\item if $A$ has the property $Q$, then so does $\fu{A}$;
\item if $A_1$ and $A_2$ have the property $Q$, then so does $A_1\oplus_{\rm alg} A_2$;
\item if $A$ has the property $Q$ and $B$ is a closed subalgebra of $A$ with finite codimension, then $B$ has the property $Q$.
\end{itemize}
Then the property $Q$ is preserved by taking Lau products in the sense of Definition~\ref{d:lauprod}.
\end{prop}

Note that two examples of such properties are: Arens regularity; and being isomorphic to a closed subalgebra of ${\mathcal B}(H)$ for some Hilbert space~$H$.

\subsection*{Acknowledgements}
The author thanks M. Nemati for pointing out the reference~\cite{JaNe_pre}, and thanks L. Molnar for helpful exchanges.
He also thanks the referee for a suggestion which improved Example~\ref{eg:different}.
\end{section}

\bibliography{nogenlauprod}
\bibliographystyle{siam}

\vfill

\noindent
Department of Mathematics and Statistics\\
Lancaster University\\
Lancaster, LA1 4YF\\
United Kingdom

\noindent
Email: {\tt y.choi.97@cantab.net}

\end{document}